\documentclass[11pt]{amsart}
\usepackage{amssymb,amsmath,amsgen,amsxtra,amsfonts} 
\usepackage{dsfont,times}

\begin{document}
%
%
%
\theoremstyle{definition}
\newtheorem{Definition}{Definition}[section]
\newtheorem*{Definitionx}{Definition}
\newtheorem{Convention}{Definition}[section]
\newtheorem{Construction}{Construction}[section]
\newtheorem{Example}[Definition]{Example}
\newtheorem{Examples}[Definition]{Examples}
\newtheorem{Remark}[Definition]{Remark}
\newtheorem{Remarks}[Definition]{Remarks}
\newtheorem{Caution}[Definition]{Caution}
\newtheorem{Conjecture}[Definition]{Conjecture}
\newtheorem{Question}[Definition]{Question}
\newtheorem*{Question*}{Question}
\newtheorem*{Acknowledgements}{Acknowledgements}
\newtheorem{Questions}[Definition]{Questions}
\theoremstyle{plain}
\newtheorem{Theorem}[Definition]{Theorem}
\newtheorem*{Theorem*}{Theorem}
\newtheorem{Proposition}[Definition]{Proposition}
\newtheorem*{Propositionx}{Proposition}
\newtheorem{Lemma}[Definition]{Lemma}
\newtheorem{Corollary}[Definition]{Corollary}
\newtheorem*{Corollary*}{Corollary}
\newtheorem{Fact}[Definition]{Fact}
\newtheorem{Facts}[Definition]{Facts}
\newtheoremstyle{voiditstyle}{3pt}{3pt}{\itshape}{\parindent}%
{\bfseries}{.}{ }{\thmnote{#3}}%
\theoremstyle{voiditstyle}
\newtheorem*{VoidItalic}{}
\newtheoremstyle{voidromstyle}{3pt}{3pt}{\rm}{\parindent}%
{\bfseries}{.}{ }{\thmnote{#3}}%
\theoremstyle{voidromstyle}
\newtheorem*{VoidRoman}{}

%
\newcommand{\prf}{\par\noindent{\sc Proof.}\quad}
\newcommand{\blowup}{\rule[-3mm]{0mm}{0mm}}
\newcommand{\cal}{\mathcal}
\newcommand{\Aff}{{\mathds{A}}}
\newcommand{\BB}{{\mathds{B}}}
\newcommand{\CC}{{\mathds{C}}}
\newcommand{\FF}{{\mathds{F}}}
\newcommand{\GG}{{\mathds{G}}}
\newcommand{\HH}{{\mathds{H}}}
\newcommand{\NN}{{\mathds{N}}}
\newcommand{\ZZ}{{\mathds{Z}}}
\newcommand{\PP}{{\mathds{P}}}
\newcommand{\QQ}{{\mathds{Q}}}
\newcommand{\RR}{{\mathds{R}}}
\newcommand{\Liea}{{\mathfrak a}}
\newcommand{\Lieb}{{\mathfrak b}}
\newcommand{\Lieg}{{\mathfrak g}}
\newcommand{\Liem}{{\mathfrak m}}
\newcommand{\ideala}{{\mathfrak a}}
\newcommand{\idealb}{{\mathfrak b}}
\newcommand{\idealg}{{\mathfrak g}}
\newcommand{\idealm}{{\mathfrak m}}
\newcommand{\idealp}{{\mathfrak p}}
\newcommand{\idealq}{{\mathfrak q}}
\newcommand{\idealI}{{\cal I}}
\newcommand{\lin}{\sim}
\newcommand{\num}{\equiv}
\newcommand{\dual}{\ast}
\newcommand{\iso}{\cong}
\newcommand{\homeo}{\approx}
\newcommand{\mm}{{\mathfrak m}}
\newcommand{\pp}{{\mathfrak p}}
\newcommand{\qq}{{\mathfrak q}}
\newcommand{\rr}{{\mathfrak r}}
\newcommand{\pP}{{\mathfrak P}}
\newcommand{\qQ}{{\mathfrak Q}}
\newcommand{\rR}{{\mathfrak R}}
%
%
\newcommand{\dq}{{``}}
\newcommand{\OO}{{\cal O}}
\newcommand{\into}{{\hookrightarrow}}
\newcommand{\onto}{{\twoheadrightarrow}}
\newcommand{\Spec}{{\rm Spec}\:}
\newcommand{\BigSpec}{{\rm\bf Spec}\:}
\newcommand{\Proj}{{\rm Proj}\:}
\newcommand{\Pic}{{\rm Pic }}
\newcommand{\Br}{{\rm Br}}
\newcommand{\NS}{{\rm NS}}
\newcommand{\chit}{\chi_{\rm top}}
\newcommand{\KanDiv}{{\cal K}}
\newcommand{\perdef}{{\stackrel{{\rm def}}{=}}}
\newcommand{\Cycl}[1]{{\ZZ/{#1}\ZZ}}
\newcommand{\Sym}{{\mathfrak S}}
\newcommand{\Xcan}{X_{{\rm can}}}
\newcommand{\Ycan}{Y_{{\rm can}}}
\newcommand{\ab}{{\rm ab}}
\newcommand{\Aut}{{\rm Aut}}
\newcommand{\Lie}{{\rm Lie}}
\newcommand{\Hom}{{\rm Hom}}
\newcommand{\Supp}{{\rm Supp}}
\newcommand{\ord}{{\rm ord}}
\newcommand{\shA}{{\mathcal{A}}}
\newcommand{\shAut}{{\underline{\rm Aut}}}
\newcommand{\shHom}{{\underline{\rm Hom}}}
\newcommand{\divisor}{{\rm div}}
\newcommand{\Alb}{{\rm Alb}}
\newcommand{\Jac}{{\rm Jac}}
\newcommand{\Ig}{{\rm Ig}}
\newcommand{\piet}{{\pi_1^{\rm \acute{e}t}}}
\newcommand{\Het}[1]{{H_{\rm \acute{e}t}^{{#1}}}}
\newcommand{\Hfl}[1]{{H_{\rm fl}^{{#1}}}}
\newcommand{\Hcris}[1]{{H_{\rm cris}^{{#1}}}}
\newcommand{\HdR}[1]{{H_{\rm dR}^{{#1}}}}
\newcommand{\hdR}[1]{{h_{\rm dR}^{{#1}}}}
\newcommand{\defin}[1]{{\bf #1}}

\title[$p$-torsion of elliptic curves]{The $\mathbf p$-torsion subgroup scheme of an elliptic curve}

\date{July 12, 2011}

\author[Christian Liedtke]{Christian Liedtke}
\address{Department of Mathematics, Stanford University, 450 Serra Mall, Stanford, CA 94305, USA}
\curraddr{}
\email{liedtke@math.stanford.edu}

\maketitle

\begin{abstract}
 Let $k$ be a field of positive characteristic $p$.

 Question: Does every twisted form of $\mu_p$ over $k$ occur as subgroup scheme
 of an elliptic curve over $k$?

 Then we show that this is true for most finite fields, for local fields and for 
 fields of characteristic $p\leq11$.
 However, it is false in general for fields of characteristic $p\geq13$, 
 which implies that there are also $p$-divisible and formal groups of 
 height one over such fields that do not arise from elliptic curves.
 It also implies that the Hasse invariant does 
 not obey the Hasse principle.

 Moreover, we also analyse twisted forms of $p$-torsion subgroup schemes
 of ordinary elliptic curves and the analogous questions for
 supersingular curves.
\end{abstract}

\setcounter{tocdepth}{1}
\maketitle
\renewcommand{\labelenumi}{(\roman{enumi})}

\section*{Introduction}

By a fundamental theorem of Raynaud, 
every finite flat and Abelian group scheme over 
some base scheme  $S$ 
can be embedded Zariski-locally into an Abelian scheme,
see \cite[Th\'eor\`eme 3.1.1]{bbm}.
%
Now, this theorem does not give a bound on the dimension
of the Abelian scheme
and a natural question would be what kind of group schemes can be realised,
say, in elliptic curves over $S$.
For example, if $S$ is some scheme of positive characteristic $p$,
and $E$ is an ordinary elliptic curve over $S$, then $\ker(F)$,
the kernel of Frobenius $F:E\to E^{(p)}$, is a finite flat and infinitesimal
group scheme over $S$, which is a twisted form of $\mu_{p,S}$.
As usual, $\mu_p$ denotes the subgroup scheme of $p$.th roots of
unity of $\GG_m$.

%

Since twisted forms of $\mu_p$ and elliptic curves are such fundamental objects
in algebraic geometry it would be natural to expect that for a broad class
of base schemes $S$ of positive characteristic $p$ every twisted form
of $\mu_p$ over $S$ can be realised in an ordinary elliptic curve over
$S$.
Being more modest, we restrict ourselves to fields and ask

\begin{Question*}[$\mathbf{A}$]
  Given a field $k$ of positive characteristic $p$, does every twisted form
  of $\mu_p$ over $k$ occur as subgroup scheme of an elliptic curve over $k$?
\end{Question*}

We note that twisted forms of $\mu_p$ over some base scheme $S$ 
are classified by the
\'etale cohomology group $\Het{1}(S, \Aut({\mu_p}))$.

Another good reason for concentrating on twisted forms of $\mu_p$ inside
elliptic curves is the following:
namely, determining
$\ker(F)$ is equivalent to knowing the Hasse invariant of
the elliptic curve.
Thus, Question (A) can be rephrased as asking whether every element of
$\Het{1}(\Spec k, \Aut({\mu_p}))$ occurs as Hasse invariant 
of an elliptic curve over $k$.

We will prove that this question has a positive answer in the following
cases:

\begin{Theorem*}
 If a field of positive characteristic $p$ is
 \begin{enumerate}
  \item a finite field with $p^n$ elements such that $p\leq17$ or $n\geq2$, or
  \item a field of the form $k((t))$, or 
  \item of characteristic $p\leq11$,
 \end{enumerate}
 then Question (A) has a positive answer for this field.
\end{Theorem*}

To get this result, 
we use Honda--Tate theory 
for (i),
a result of Igusa on Hasse invariants of 
semi-universal deformations of supersingular elliptic curves
to obtain (ii),
and (iii) follows from straight forward computations with
Hasse invariants.
Moreover, we even get the following analogue of Raynaud's theorem:

\begin{Theorem*}
 Let $S$ be a scheme of positive characteristic $p\leq11$.
 Then, every twisted form of $\mu_{p,S}$ over $S$ can be
 embedded Zariski-locally into an elliptic curve.
\end{Theorem*}

Conversely, we will prove that these results are rather sharp:

\begin{Theorem*}
  If a  field of positive characteristic $p$ is
  \begin{enumerate}
  \item a prime field $\FF_p$ with $p\geq 19$, or
  \item of the form $k(t)$ with $p\geq13$ 
  \end{enumerate}
  then Question (A) has a negative answer for this field.
\end{Theorem*}

The rather odd appearance of the prime $13$ in (ii) is closely 
related to the fact
that the Igusa curves $\Ig(p)^{\rm ord}$ are non-rational curves if and only if
$p\geq13$, and this turns out to be an obstruction to realising twisted
forms of $\mu_p$ over $k(t)$.

Moreover, we will see 
that we cannot improve this theorem by allowing twisted
forms of $\mu_p$ on nonproper or singular curves carrying 
group structures.

A direct consequence is

\begin{Corollary*}
 For every $p\geq13$ (resp. $p\geq19$)
 there exists a global (resp. finite) field $K$ 
 of characteristic $p$ and
 \begin{enumerate}
  \item a twisted form of $\mu_p$ over $K$,
  \item a $p$-divisible group of height one over $K$, and 
  \item a formal group of height one over $K$, 
 \end{enumerate}
 such that neither of them arises from an elliptic curve over $K$.
\end{Corollary*}

From class field theory it follows that twisted forms
of $\mu_p$ over global fields obey a Hasse principle.
However, the theorems above imply

\begin{Corollary*}
 For every $p\geq13$ there exists a global field $K$ of 
 characteristic $p$
 such that the Hasse principle does not hold for 
 Hasse invariants of elliptic curves over $K$.
\end{Corollary*}
\medskip

Moreover, the $p$-torsion subgroup scheme $E[p]$ of an elliptic curve sits 
inside an extension
\begin{equation*}
 \tag{\ref{connectedetale}}
 0\,\to\,\ker(F)\,\to\,E[p]\,\to\,\ker(V)\,\to\,0
\end{equation*}
(here, $V:E^{(p)}\to E$ denotes Verschiebung), 
and is endowed with a bilinear, alternating
and nondegenerate pairing, the Weil pairing.
Due to this pairing, kernel and cokernel of (\ref{connectedetale}) become
Cartier dual group schemes.
Hence we may ask

\begin{Question*}[$\mathbf{B}$]
 Given a field $k$ of positive characteristic $p$,
 does every twisted form of $\mu_p\oplus(\ZZ/p\ZZ)$ over $k$ that is endowed with a
 bilinear, alternating and nondegenerate pairing occur as subgroup scheme
 of an elliptic curve over $k$?
\end{Question*}

For perfect fields, the sequence (\ref{connectedetale}) is split from which it follows
that Question (A) and Question (B) are equivalent for these fields.
Also, if already Question (A) has a negative answer for a field,
Question (B) cannot have a positive answer.
In view of our positive results we can only hope for local fields and 
for fields of characteristic $p\leq11$ to have a positive answer to Question (B),
and in fact:

\begin{Theorem*}
 If a field of positive characteristic $p$ is 
 \begin{enumerate}
  \item a field of the form $k((t))$ where $k$ is perfect, or 
  \item if $p\leq11$,
 \end{enumerate}
 then Question (B) has a positive answer for this field.
\end{Theorem*}

Finally, we answer the analogous questions for supersingular curves:
if $E$ is supersingular then
$\ker(F)$ is a twisted form of $\alpha_p$,
and $E[p]$ is a twisted form of $M_2$.
These latter group schemes denote the 
the unique non-split extension of $\alpha_p$ by itself
that is autodual, cf.
\cite[Section II.15.5]{oort}.
On the other hand, we will see that neither $\alpha_p$ nor
$M_2$ possess twisted forms over fields, whence
$\ker(F)\iso\alpha_p$ and $E[p]\iso M_2$ for every supersingular
elliptic curve.
In particular, Questions (A) and (B) also hold true for supersingular
elliptic curves, albeit for trivial reasons.
\medskip

This article is organised as follows:

In Section \ref{sec:generalities} we recall a couple of facts
from \cite{kama} and \cite{ls}
about twisted forms of $\mu_p$, Hasse invariants and $p$-torsion
subgroup schemes of ordinary elliptic curves.

In Section \ref{sec:finite} we use Weil's results
on elliptic curves over finite fields and Honda--Tate theory
\cite{tate} to answer Questions (A) and (B) for finite
fields.

In Section \ref{sec:local} we use universal formal deformations
of supersingular elliptic curves and a result of Igusa on the
Hasse invariants of such deformations to get a positive
answer to Question (A) for all fields of the form $k((t))$.
In particular, this gives a positive answer to Question (A) 
for local fields.
Finally, we prove a Hasse principle for twisted forms over
global fields using class field theory.

In Section \ref{sec:igusa} we first answer Question (A) positively
for all fields of characteristic $p\leq11$ using explicit
computations with Hasse invariants.
Then we prove that Question (A) has a negative answer for 
the function field of $\PP^1$ in characteristic $p\geq13$ and 
for function fields of elliptic curves in characteristic 
$p\geq17$.
This result is closely related to Igusa's result \cite{igusa} 
that the Igusa curves $\Ig(p)$ are not rational for $p\geq13$
and of general type for $p\geq17$.

In Section \ref{sec:ptorsion}
we show that Question (B) also has a positive answer for
local fields and for fields of characteristic $p\leq11$.

Finally, in Section \ref{sec:supersingular} we discuss
kernel of Frobenius and the $p$-torsion subgroup scheme
of supersingular elliptic curves over fields.

\begin{Acknowledgements}
 I thank Gerhard~Frey, Frans~Oort and the referee for comments and 
 pointing out some of the references to me. 
\end{Acknowledgements}

\section{Generalities}
\label{sec:generalities}

In this section we recall a couple of facts concerning the kernel of Frobenius, Hasse invariants
and $p$-torsion subgroup schemes of ordinary elliptic curves.
For generalities about torsion subgroup schemes of elliptic curves we refer to \cite[Chapter 8.7]{kama}.
For supersingular elliptic curves, see Section \ref{sec:supersingular} below.

Let $S$ be a base scheme of characteristic $p>0$ and $E$ be an elliptic curve over $S$.
Then $\ker(F)$, the kernel of Frobenius $F:E\to E^{(p)}$, is a finite and flat group scheme of length
$p$ over $S$.
The {$p$-Lie algebra} of $\ker(F)$ is a projective $\OO_S$-module of rank $1$ and thus coincides
with the $p$-Lie algebra of $E$.

We will now assume that $S=\Spec T$ is affine with $\Pic(T)=0$, e.g., $T$ could be a field, a local
ring or a polynomial ring over a field.
In this case a sheaf of projective $\OO_S$-modules of rank $1$ admits a global basis $t\in\OO_S$
identifying this module with $T$.
Moreover, giving $S$ the structure of a sheaf of $p$-Lie algebras is equivalent to
prescribing the action of the $p$-operator on this basis, i.e.,
$t^{[p]} = \lambda t$ for some $\lambda\in T$.
Note however, that we have chosen a basis for $\OO_S$, which makes $\lambda$ only unique
up to multiplication by elements of $T^{\times(p-1)}$, confer
\cite[Chapter V, Sections 7 and 8]{jacobson}.
Under our assumptions on $S$,
it follows from the Tate--Oort classification \cite{to} of group schemes of prime order 
that the $p$-Lie algebra of $\ker(F)$ determines $\ker(F)$ uniquely.

Let us now assume that $E$ is an {ordinary elliptic curve} over $S=\Spec T$, 
which is equivalent to $\lambda\in T^\times$ in our setup.
Then $\ker(F)$ is a {twisted form of $\mu_p$} over $S$.
Since the automorphism group scheme $\shAut(\mu_p)$ of $\mu_p$ is isomorphic 
to $\mu_{p-1}$, twisted forms of $\mu_p$ are classified by
$\Het{1}(S,\mu_{p-1})$.
Using the Kummer sequence $0\to\mu_{p-1}\to\GG_m\to\GG_m\to0$ and our assumption
$\Pic(T)=0$ we deduce that the coboundary map in cohomology induces 
an isomorphism of groups
$$
  T^\times/T^{\times(p-1)}\,\to\,\Het{1}(S,\mu_{p-1})\,.
$$
Here we see directly that the $p$-Lie algebra determines the twisted form of
$\mu_p$ uniquely.
See also \cite[Section 1]{ls} for the preceding discussion.

If we drop the assumption on $\Pic(T)$ for a moment, then
a twisted form of $\mu_p$ over $S$ determines - via the long exact
sequence in cohomology associated to the Kummer sequence - an element
in $H^1(S,\GG_m)$, i.e., an invertible sheaf.
After trivialising this sheaf by
passing to an affine Zariski-open cover $\{U_\alpha\}_{\alpha\in I}$ 
of $S$, the restriction of the twisted form of $\mu_p$ to $U_\alpha$ is
given by a class in $T_\alpha/T_\alpha^{\times(p-1)}$, where 
$T_\alpha=H^0(U_\alpha,\OO_{U_\alpha})$.

We recall from \cite[Section 12.4]{kama} that the {\em Hasse invariant}
of an elliptic curve is defined to be the linear mapping induced by
the Verschiebung $V$ on $p$-Lie algebras:
$h(E)=\Lie(V):\Lie(E^{(p)})\to\Lie(E)$.
Using the identification $\Lie(E^{(p)})=\Lie(E)^{\otimes p}$, we may regard
the Hasse invariant as an element of the projective $\OO_S$-module
$\Lie(E)^{\otimes(1-p)}$, which is of rank one.
Choosing a basis $u\in\Lie(E)$, we can identify $h(E)$ with 
$\lambda' u^{\otimes(1-p)}$
for some $\lambda'\in T$ which is unique up to $(p-1)$.st powers.
As carried out in \cite[Section 3]{ls}
the Hasse invariant determines the $p$-Lie algebra of $E$ up to
isomorphism and conversely.
More precisely, we have $\lambda'=\lambda^{-1}$ by
\cite[Proposition 3.2]{ls}.
In particular, Question (A) of the introduction is equivalent to 

\begin{Question*}[$\mathbf{A'}$]
 Given a field $k$ of positive characteristic $p$, does every element of 
 $k^\times/k^{\times(p-1)}\iso\Het{1}(\Spec k,\mu_{p-1})$ occur as Hasse invariant
 of an ordinary elliptic curve over $k$?
\end{Question*}

To describe the $p$-torsion subgroup scheme $E[p]$ of an ordinary elliptic curve 
$E$ over $S$, let us recall the setup of \cite[Section 3]{ls}:
this group scheme is a twisted form of $G:=\mu_p\oplus(\ZZ/p\ZZ)$.
On $E[p]$ there exists a bilinear, alternating and nondegenerate pairing,
the {\em Weil pairing}.
This pairing makes the identity component $\ker(F)$ of $E[p]$ the
Cartier dual of the \'etale quotient $\ker(V)$
\begin{equation}
 \label{connectedetale}
  0\,\to\,\ker(F)\,\to\,E[p]\,\to\,\ker(V)\,\to\,0\,.
\end{equation}
In particular, not every twisted form of $G$ can occur as $p$-torsion subgroup
scheme of an ordinary elliptic curve.
However, on $G$ we may define a bilinear, alternating and nondegenerate
pairing by
$$
\Phi\,:\,G\times G\,\to\,\mu_p\,\mbox{ \quad } ((\mu,i),\,(\nu,j))\,\mapsto\,\mu^j/\nu^i\,,
$$
and denote by $\shA:=\shAut(G,\Phi)$ the automorphism scheme of those automorphisms of $G$
that respect $\Phi$.
Then $E[p]$ defines a cohomology class in $\Hfl{1}(S,\shA)$.
We can now reformulate Question (B) of the introduction:

\begin{Question*}[$\mathbf{B'}$]
 Given a field $k$ of positive characteristic $p$, can every cohomology class in
 $\Hfl{1}(\Spec k,\shA)$ be realised by the $p$-torsion subgroup scheme of an 
 ordinary elliptic curve over $k$?
\end{Question*}

In \cite[Section 2]{ls} we described $\shA$ explicitly: it sits inside a split short
exact sequence
\begin{equation}
 \label{autsequence}
 1\,\to\,\mu_{p-1}\,\to\,\shA\,\to\,\mu_p\,\to\,1\,.
\end{equation}
It turns out that $\shA$ is not commutative for $p\geq5$ and then non-Abelian group 
cohomology is needed to describe twisted forms of $(G,\Phi)$.
Taking cohomology in (\ref{autsequence}) induces a surjective homomorphism of
pointed sets
\begin{equation}
 \label{aut cohomology}
  \Hfl{1}(S,\,\shA)\,\to\,\Hfl{1}(S,\,\mu_{p-1})\,\to\,1
\end{equation}
mapping the class of $E[p]$ to the class of $\ker(F)$.
Since $\shA$ is non-Abelian for $p\geq5$, the description of the fibres of
(\ref{aut cohomology}) is a little bit tricky.
Roughly speaking, these fibres describe extension classes like the 
short exact sequence (\ref{connectedetale}).
We refer to \cite[Theorem 3.4]{ls} for details.

The sequence (\ref{connectedetale}) is just the connected-\'etale exact sequence
for $E[p]$ and it is well known that such sequences split over perfect fields.
In our case this can also be seen from (\ref{aut cohomology}), since
its kernel $\Hfl{1}(\Spec k,\mu_p)$ is trivial for perfect fields.
More generally, (\ref{connectedetale}) splits if and only if 
$j(E)\in k^p$ by \cite[Proposition 3.3]{ls}.
In any case we note

\begin{Remark}
 \label{perfect remark}
 For a perfect field, Questions (A) and (B) are equivalent.
\end{Remark}

\section{Finite fields}
\label{sec:finite}

In this section we answer Questions (A) and (B) for finite fields.
After classifying twisted forms over finite fields we use Weil's results
about counting points of elliptic curves over finite fields and Honda--Tate
theory \cite{tate} to decide which twisted forms of $\mu_p$ over a finite
field can be realised as subgroup schemes of elliptic curves
over this finite field.

\begin{Proposition}
 \label{classification finite}
 Let $k=\FF_q$ be the finite field with $q=p^n$ elements. 
 Then there exists a bijection of sets and an isomorphism $\varphi$ of Abelian
 groups
 $$
   \{\mbox{ twisted forms of $\mu_p$ over $k$ }\} \,\to\,
   k^\times/k^{\times(p-1)}\,\stackrel{\varphi}{\to}\,\FF_p^\times\,.
 $$
 In particular, there are precisely $(p-1)$ twisted forms of $\mu_p$ over $k$.
\end{Proposition}

\prf
We have seen in Section \ref{sec:generalities} that
twisted forms of $\mu_p$ over $k$ correspond bijectively to one-dimensional
$p$-Lie algebras over $k$ whose $p$-operator is non-trivial.
There we have also seen that these $p$-Lie algebras correspond bijectively
to the set $k^\times/k^{\times(p-1)}$.
Finally, we leave it to the reader to check that
\begin{equation}
 \label{varphi def}
 \begin{array}{cccccc}
    \varphi &:& k^\times/k^{\times(p-1)}&\to&\FF_p^\times\\
    &&x&\mapsto&x^m&\mbox{ \quad with \quad } m=\frac{q-1}{p-1}
  \end{array}
\end{equation}
defines an isomorphism of Abelian groups.
\qed

\begin{Theorem}
 \label{positive finite}
 Let $k=\FF_q$ be the finite field with $q=p^n$ elements and assume that
 $p\leq17$ or $n\geq2$.
 Then every twisted form of $\mu_p$ over $k$ occurs as subgroup scheme of
 an elliptic curve over $k$.
 In particular, Questions (A) and (B) have a positive answer for $k$.
\end{Theorem}

\begin{Theorem}
 \label{negative finite}
 Let $\FF_p$ be a prime field with $p\geq19$ elements.
 Then the set of twisted forms that occur as subgroup schemes of elliptic
 curves over $\FF_p$ corresponds via the map $\varphi$ of 
 Proposition \ref{classification finite} to the set
 $$
  \{\,\, [\beta]\in\FF_p^\times\,|\, \beta\in\ZZ-\{0\},\,\beta^2<4p\,\,\}
  \,\subset\,\FF_p^\times\,,
 $$
 which is a proper subset.
 In particular, Questions (A) and (B) have a negative answer for $\FF_p$.
\end{Theorem}

{\sc Proof (of both theorems)}. 
Since finite fields are perfect, Questions (A) and (B) are equivalent
by Remark \ref{perfect remark}.
By Proposition \ref{classification finite} there are no twisted forms of
$\mu_p$ in characteristic $p=2$ and we may thus assume $p\geq3$.

Let $E$ be an elliptic curve over the finite field $k=\FF_q$ with
$q=p^n$ elements.
Then the characteristic polynomial on the first $\ell$-adic cohomology
group is of the form
$x^2-\beta x+q$,
where $\beta$ is an integer satisfying $\beta^2<4q$, confer 
\cite[Chapter V, Section 2]{aec}.
By loc.cit. we also know that the number of $k$-rational points on
$E$ is equal to 
$$
\# E(k) \,=\, 1\,-\,\beta\,+\,q\,.
$$
Since $p\geq3$ we may assume that $E$ is given by an equation of the
form $y^2=f(x)$.
Denote by $A_{p^m}$ the coefficient of $x^{p^m-1}$ in $f(x)^{(p^m-1)/2}$.
Then $A_p$ is the Hasse invariant of $E$ and its class in 
$k^\times/k^{\times(p-1)}$ corresponds to the twisted form $\ker(F)$
of $\mu_p$ as in Proposition \ref{classification finite} and as explained
in Section \ref{sec:generalities}.
From the proof of \cite[Theorem V.4.1]{aec} we get
$$
A_q\,=\,A_p^{1+p+...+p^{n-1}}\,=\,A_p{}^{\frac{q-1}{p-1}}\,=\,\varphi(A_p),
$$
where $\varphi$ is as in (\ref{varphi def}).
In particular, $A_q$ lies in $\FF_p$.
From loc.cit. we also get 
$$
A_q\,\equiv\,1\,-\,\#E(k)\,\mod\,p\,.
$$
Putting these results together, we infer
$$
\beta\,=\,1\,-\,\#E(k)\,+\,q\,\equiv\,A_q\,=\,\varphi(A_p)\,\mod\,p,
$$
i.e., $\beta\mod p$ as an element of $\FF_p$ determines $\ker(F)$ via
the correspondence of Proposition \ref{classification finite}.

Theorem \ref{negative finite} now follows:
in fact, $\beta$ is an integer fulfilling $\beta^2<4q$, and if
$q=p\geq19$ there are strictly less than $p-1$ possibilities for $\beta$.
On the other hand, there are precisely $p-1$ twisted forms of $\mu_p$
over $\FF_p$ and hence not every one of them can be realised as
kernel of Frobenius on an elliptic curve.

Conversely, for $h\in\FF_p^\times$ and $q=p^n$ with $p\leq17$ or $n\geq2$
there exists an integer $\beta$ with $\beta^2<4q$ and $\beta\equiv h\mod p$.
Now, Honda--Tate theory tells us that there exists an elliptic curve over
$k=\FF_q$ with characteristic polynomial of Frobenius equal to
$x^2-\beta x+q$, confer \cite{tate}.
This proves Theorem \ref{positive finite}.
\qed

\begin{Remark}
 \label{mup ok}
 The class $[\beta]=1\in\FF_p^\times$ corresponds to $\mu_p$ and thus
 Honda--Tate theory shows moreover that there always exists an elliptic
 curve $E$ over $\FF_p$ with $\ker(F)\iso\mu_p$ and thus
 $E[p]\iso\mu_p\oplus(\ZZ/p\ZZ)$.
\end{Remark}

\begin{Remark}
 From the description in Theorem \ref{negative finite} 
 of those twisted forms of 
 $\mu_p$ that are realisable as subgroup schemes of elliptic curves 
 it is not difficult to see that this subset usually does not form a
 subgroup.
\end{Remark}

\section{Local fields}
\label{sec:local}

In this section we show that Question (A) has a positive answer for all fields
of the form $k((t))$ using deformation theory of elliptic curves 
and a result of Igusa.
Also, we prove a Hasse principle for twisted forms of $\mu_p$
and classify these forms over local fields.

\begin{Theorem}
 \label{local positive}
 Let $k$ be a field of positive characteristic $p$.
 Then every twisted form of $\mu_p$ over $k((t))$ occurs as subgroup scheme
 of an elliptic curve over $k((t))$. 
 In particular, Question (A) has a positive answer for $k((t))$.
\end{Theorem}

\prf
Choose a supersingular elliptic curve over $\FF_p$ and let
${\cal E}\to\FF_p[[t]]$ be its universal formal deformation.
By a theorem of Igusa (\cite[Corollary 12.4.4.]{kama}), the Hasse
invariant of $\cal E$ has simple zeroes, i.e., the Hasse invariant
can be lifted to an element of $\FF_p[[t]]$ of the form
$$
c\cdot t\cdot t^{n(p-1)}\cdot(1+tg(t))
$$
for some $c\in\FF_p^\times$, some $n\geq0$ and $g(t)\in\FF_p[[t]]$.
Using the fact that the equation $x^{p-1}-1$ has $p-1$ different zeroes 
in $\FF_p$ and using Hensel's lemma, it is easy to see that
this Hasse invariant is congruent to $c\cdot t$ modulo
$\FF_p((t))^{\times(p-1)}$.

Base changing from $\FF_p$ to $k$ and using then base changes of the
form $t\mapsto a\cdot t^{m}$ we can realise every 
class in $k((t))^\times/k((t))^{\times(p-1)}$ as
Hasse invariant of an elliptic curve over $k((t))$.
\qed\medskip

Now, let $K$ be a global field of positive characteristic $p$, i.e.,
a finite extension of $\FF_p(t)$.
For a place $v$ of $K$ we denote by $K_v$ the completion 
with respect to $v$.
The next result shows that twisted forms of $\mu_p$ over
global fields obey a Hasse principle.

\begin{Proposition}[Hasse principle]
 \label{hasse principle}
 Let $K$ be global field of positive characteristic $p$.
 Then the natural homomorphism 
 $$
  \Het{1}(\Spec K,\,\mu_{p-1})\,\to\,\prod_v \Het{1}(\Spec K_v,\,\mu_{p-1})\,,
 $$
 where the product is taken over all places of $K$, is injective.
 In particular, a twisted form of $\mu_p$ over $K$ is trivial if and only
 if it is trivial over every completion.
\end{Proposition}

\prf
The total space of a $\mu_{p-1}$-torsor over $\Spec K$ is the spectrum of
an Artin algebra over $K$, hence a product of fields $L\supset K$.
If this torsor is non-trivial then $L/K$ is a non-trivial extension
and due to the $\mu_{p-1}$-action this extension is a Galois extension
with Abelian Galois group.
By class field theory there exists a place $v$ of $K$ such that the
induced extension $L_v/K_v$ on completions is non-trivial.
This implies that the induced $\mu_{p-1}$-torsor over $K_v$ is non-trivial
and our injectivity statement follows.
\qed\medskip

\begin{Remark}
 Let $K=k((t))$ for some field $k$ of positive characteristic $p$.
 Using Hensel's lemma it is easy to see that the valuation
 $\nu:K^\times\to\ZZ$ induces a short exact sequence
 $$
  0\,\to\,\Het{1}(\Spec k,\mu_{p-1})\,\to\,\Het{1}(\Spec K,\mu_{p-1})
  \,\stackrel{\nu}{\to}\,\ZZ/(p-1)\ZZ\,\to\,0\,,
 $$
 which can be split using the uniformiser $t\in K$.
 In particular, if $K$ is a local field, i.e., if $k$ is a 
 finite field $\FF_q$, then applying
 Proposition \ref{classification finite} 
 we obtain an isomorphism
 $$
    \Het{1}(\Spec \FF_q((t)),\,\mu_{p-1})\,\iso\,(\ZZ/(p-1)\ZZ)^2\,.
 $$
 In particular, there are only finitely many twisted forms of $\mu_p$
 over local fields.
\end{Remark}

\section{Igusa curves}
\label{sec:igusa}

In this section we will first prove that Question (A) has a positive answer
not only for all fields, but even
Zariski-locally for arbitrary schemes of positive characteristic $p\leq11$.
Then, however, we will see that Question (A) has a negative answer
for the fields $k(t)$ in characteristic $p\geq13$ and for 
function fields of elliptic curves in characteristic $p\geq17$.
This result is closely related to the geometry of the Igusa curves.

\begin{Theorem}
  \label{easy main result}
  Let $S$ be a scheme of positive characteristic $p\leq11$.
  Then, every twisted form of $\mu_{p,S}$ over $S$ admits Zariski-locally
  an embedding into an elliptic curve.
  That is, there exists a Zariski-open cover $\{U_\alpha\}_{\alpha\in I}$
  of $S$ and elliptic curves $E_\alpha$ over $U_\alpha$, such that
  the restriction $\widetilde{\mu}_{p,S}|_{U_\alpha}$ embeds into $E_\alpha$ for
  all $\alpha\in I$
\end{Theorem}

\prf
A twisted form $\widetilde{\mu}_{p,S}$ of $\mu_{p,S}$ corresponds to 
a class in $\Het{1}(S,\mu_{p-1})$, see Section \ref{sec:generalities}.
The long exact sequence in cohomology associated to the
Kummer sequence $0\to\mu_{p-1}\to\GG_m\to\GG_m\to0$ yields
$$
\begin{array}{cccccccc}
 0 & \to & S^\times/S^{\times(p-1)}
 &\stackrel{\delta}{\to}& 
 \Het{1}(S,\mu_{p-1}) &\to&H^1(S,\GG_m)&...\\
 &&&&[\widetilde{\mu}_{p,S}]&\mapsto&{\cal L}
\end{array}
$$
We choose an open affine cover $\{V_\beta\}_{\beta\in J}$ 
of $S$
trivialising the invertible sheaf $\cal L$.
Then the restriction of $\widetilde{\mu}_{p,S}$ to
$V_\beta$ corresponds to a class in
$T_\beta^\times/T_\beta^{\times(p-1)}$, where
$T_\beta=H^0(V_\beta,\OO_{V_\beta})$, see
also Section \ref{sec:generalities}.

If $p\geq5$ then every elliptic curve over a field 
can be given by a Weierstra\ss\ equation
of the form $y^2=x^3+ax+b$.
Depending on $p$ we obtain the following Hasse invariants:
$$
\begin{array}{c|ccc}
 p & 5 & 7 & 11 \\
 \hline
 \mbox{Hasse invariant} & 2a & 3b & 9ab
\end{array}
$$
Hence if $5\leq p\leq11$ we can easily find for every class 
$[h_\beta]\in T_\beta^\times/T_\beta^{\times(p-1)}$
and every closed point $Q\in V_\beta$
an affine Zariski-open neighbourhood $V_Q$ of $Q$ in $V_\beta$
and an elliptic curve $E_Q$ over $V_Q$
with Hasse invariant $[h_\beta]\in T_Q^\times/T_Q^{\times(p-1)}$,
where $T_Q=H^0(V_Q,\OO_{V_Q})$.
Thus, $E_Q$ contains 
$\tilde{\mu}_{p,S}|_{V_Q}$ as subgroup scheme.
In particular, after possibly refining the cover $V_\beta$,
we obtain a Zariski-open cover $\{U_\alpha\}$
of $S$ and elliptic curves over $U_\alpha$ containing 
$\widetilde{\mu}_{p,S}|_{U_\alpha}$ as subgroup scheme.

We leave the cases $p=2$ and $p=3$ to the reader.
\qed\medskip

\begin{Corollary}
 \label{cor easy main result}
 Question (A) has a positive answer for fields of characteristic 
 $2\leq p\leq11$.
\end{Corollary}

\begin{Remark}
 \label{refined easy main result}
 Coming back to fields, one can ask what freedom one has choosing the 
 elliptic curve containing a given twisted form of $\mu_p$ as subgroup scheme.
 Using the automorphism groups of (special) elliptic curves to
 twist a given curve it follows from Lemma \ref{twisting} below that
 \begin{enumerate}
  \item [$p=2$] There are no twisted forms of $\mu_2$ over fields and every
      ordinary elliptic curve contains $\mu_2$ as subgroup scheme.
  \item [$p=3$] Given a twisted form $\widetilde{\mu}_3$ of $\mu_3$ and an 
     arbitrary ordinary elliptic curve $E$ there exists a quadratic twist of $E$
     containing $\widetilde{\mu}_3$ as subgroup scheme.
  \item [$p=5$] Every twisted form of $\mu_5$ can be realised as subgroup scheme
     of an elliptic curve with $j=1728$.
  \item [$p=7$] Every twisted form of $\mu_7$ can be realised as subgroup scheme
     of an elliptic curve with $j=0$.
 \end{enumerate}
 In characteristic $p=11$ there usually does not exist one single elliptic curve
 such that all twisted forms of $\mu_p$ occur as subgroup schemes of twists of
 this particular curve.
\end{Remark}

We now come to the main result of this section, namely that in characteristic
$p\geq13$ Question (A) has a negative answer in general.
As the proof will show this is closely related to the Igusa curves not
being rational if $p\geq13$.
We recall from \cite[Chapter 12.3]{kama} that the 
{\em Igusa moduli problem} 
classifies ordinary elliptic curves ${\cal E}$ over base schemes
$S$ of positive characteristic $p$
such that the Frobenius pullback ${\cal E}^{(p)}$ contains an $S$-rational
$p$-division point.
For $p\geq3$ this moduli problem is representable by a smooth affine curve
$\Ig(p)^{\rm ord}$ over $\FF_p$, whose smooth compactification is denoted by
$\overline{\Ig(p)^{\rm ord}}$.
The geometry of these curves has been analysed in \cite{igusa}, and in particular
their genera have been determined there:

\begin{Proposition}[Igusa]
 \label{igusa proposition}
 Let $\overline{\Ig(p)^{\rm ord}}$ be the smooth compactification of the
 Igusa curve in positive characteristic $p\geq3$.
 Then this curve is
 \begin{enumerate}
  \item rational if $p\leq11$, is
  \item elliptic if $p=13$, and is
  \item of general type if $p\geq17$. \qed
 \end{enumerate}
\end{Proposition}

Next, we determine the effect that
twisting an elliptic curve has on its Hasse invariant.

\begin{Lemma}
 \label{twisting}
  Let $E$ be an ordinary elliptic curve over a field $k$ of characteristic
  $p\geq3$ and let $[h]\in k^\times/k^{\times(p-1)}$ be its Hasse invariant.
  \begin{enumerate}
   \item If $E^D$ is the quadratic twist of $E$ with respect to 
      $D\in k^\times/k^{\times2}$
      then the Hasse invariant of $E^D$ is $[hD^{(p-1)/2}]$.
   \item If $j(E)=0$ (and thus $p\equiv1\mod3$) and $E^D$ is the sextic twist of $E$
      with respect to $D\in k^\times/k^{\times6}$
      then the Hasse invariant of $E^D$ is $[hD^{(p-1)/6}]$.
   \item If $j(E)=1728$ (and thus $p\equiv1\mod4$) and $E^D$ is the quartic twist of $E$
      with respect to $D\in k^\times/k^{\times4}$
      then the Hasse invariant of $E^D$ is $[hD^{(p-1)/4}]$.
  \end{enumerate}
\end{Lemma}

\prf
In characteristic $p\geq5$ this can be seen from the explicit computation of
twists in \cite[Chapter X, Proposition 5.4]{aec} together with the description
of the Hasse invariant as the coefficient of $x^{p-1}$ in $f(x)^{(p-1)/2}$ if
the elliptic curve is given by $y^2=f(x)$, see \cite[Chapter V, Theorem 4.1]{aec}.
We leave the case $p=3$ to the reader.
\qed

\begin{Theorem}
 \label{main result}
 Let $k$ be a field of characteristic $p\geq13$.
 Then there exists a twisted form of $\mu_p$ over $k(t)$ that does not occur
 as subgroup scheme of an elliptic curve over $k(t)$.
 In particular, Question (A) has a negative answer for $k(t)$.
\end{Theorem}

\prf
Consider the field extension $k(t)\subset L:=k(t)[\sqrt[p-1]{t}]$, which
is Galois with cyclic Galois group of order $p-1$.
Hence $\Spec L\to \Spec k(t)$ is a $\mu_{p-1}$-torsor
and we denote by $\widetilde{\mu}_p$ the twisted form of $\mu_p$ which arises
by twisting $\mu_p$ with this torsor.
Note that $L=k(u)$ for $u:=\sqrt[p-1]{t}$, i.e., $L$ is the function field
of $\PP^1_k$. 

By way of contradiction, we assume that there exists an elliptic curve 
$E$ over $k(t)$ containing $\widetilde{\mu}_p$ as subgroup scheme.

Suppose first that $j(E)\in k$.
Then there exists a twist $E'$ of $E$ which is already defined over $k$.
The corresponding twist $\widetilde{\mu}_p'$ is then also defined over $k$.
If $j(E)\not\in\{0,1728\}$ then $E'$ is a quadratic twist of
$E$, see \cite[Chapter X, Proposition 5.4]{aec}.
Now, $\widetilde{\mu}_p$ corresponds to the class 
$[t]\in k(t)^\times/k(t)^{\times(p-1)}$.
Since $\widetilde{\mu}_p'$ is obtained by a quadratic twist 
(still assuming $j(E)\not\in\{0,1728\}$ for the moment), 
say with twisting parameter $D$, the twisted form
of $\widetilde{\mu}_p'$ corresponds to the class $[tD^{(p-1)/2}]$
by Lemma \ref{twisting}.
Using $k[t]$ with its standard valuation $v$ such that $v(t)=1$,
we see that for $p\geq5$ and for all $D\in k(t)^\times$ we have
$v(tD^{(p-1)/2})\neq0$.
In particular, the class of $\widetilde{\mu}_p'$ in
$k(t)^\times/k(t)^{\times(p-1)}$ cannot be represented
by an element of $k$, which implies that $\widetilde{\mu}_p'$
cannot be defined over $k$.
If $j(E)=0$ or $j(E)=1728$ we have to consider also sextic and quartic
twists, but again, for $p\geq11$, no sextic or quartic twist of
$\widetilde{\mu}_p$ can be defined over $k$.
We conclude that $j(E)\not\in k$, i.e., $j(E)$ is transcendental
over $k$.

Since $\widetilde{\mu}_p$ becomes isomorphic to $\mu_p$ over $L$,
we see that $\mu_p$ is a subgroup scheme of $E_L:=E\times_{\Spec k(t)}\Spec L$.
This implies that $E_L^{(p)}$ contains an $L$-rational $p$-division
point and we obtain a classifying morphism 
$\varphi:\Spec L\to \Ig(p)^{\rm ord}$ to the Igusa curve which
yields a morphism
$\varphi_k:\Spec L\to \Ig(p)^{\rm ord}\times_{\Spec \FF_p}\Spec k$.
There is a dominant morphism $\Ig(p)^{\rm ord}\to\PP^1$ 
induced by the $j$-invariant and using that
$j(E)$ is transcendental over $k$ we infer that $\varphi_k$
is dominant.
Since $L$ is the function field of $\PP^1_k$, the existence
of $\varphi_k$ implies that $\overline{\Ig(p)^{\rm ord}}$ 
is a rational curve.
This contradicts Proposition \ref{igusa proposition} and we conclude
that the curve $E$ we started with does not exist.
\qed\medskip

This result together with Theorem \ref{negative finite}
implies that the analogue of Question (A) for formal
and $p$-divisible groups of height
one also has a negative answer in general:

\begin{Corollary}
 For every $p\geq13$ (resp. $p\geq19$)
 there exists a global (resp. finite) field $K$ 
 of characteristic $p$ and
 \begin{enumerate}
  \item a twisted form of $\mu_p$ over $K$,
  \item a $p$-divisible group of height one over $K$, and 
  \item a formal group of height one over $K$, 
 \end{enumerate}
 such that neither of them arises from an elliptic curve over $K$.
\end{Corollary}

\prf
By Theorem \ref{main result} 
(resp. Theorem \ref{negative finite})
there exists a global (resp. finite)
field $K$ and a $\mu_{p-1}$-torsor $\Spec L\to\Spec K$
such that the associated twisted form
$\widetilde{\mu}_p$ 
of $\mu_p$ cannot be realised as subgroup scheme
of an elliptic curve over $K$.

We may extend the $\mu_{p-1}$-action of $\mu_p$ to an 
action on $\mu_{p^\infty}$.
Then, twisting $\mu_{p^\infty}$ with this torsor we obtain
a $p$-divisible group of height one 
that cannot arise as $E[p^{\infty}]$
of an elliptic curve $E$ over $K$.

The formal group associated to an elliptic curve determines
its Hasse invariant and thus $\ker(F)$, see
\cite[Chapter 12.4]{kama} or Section \ref{sec:generalities}.
Let $E$ be an elliptic curve over $K$ with 
$\ker(F)\iso\mu_p$, which is always possible
by Remark \ref{mup ok}.
The associated formal group $\widehat{E}$
is of height one and thus
$\Aut(\widehat{E})\iso\ZZ_p^*$
by \cite[Proposition 9]{lazard}.
We note that $\ZZ_p^*$ naturally 
contains $\mu_{p-1}$.
Twisting $\widehat{E}$ with respect to the
$\mu_{p-1}$-torsor $\Spec L\to\Spec K$
we obtain a formal group of height one
which cannot arise from an elliptic curve
over $K$, since
it would correspond to an elliptic curve
having $\widetilde{\mu}_p$ as kernel of Frobenius.
\qed\medskip

Actually, there are many more fields than those
of Theorem \ref{main result}
that provide counter-examples to Question (A).
The following result yields infinitely many examples
for every $p\geq17$:

\begin{Theorem}
 \label{main result elliptic}
 Let $k$ be a field of characteristic $p\geq17$ and $E$ be an elliptic
 curve over $k$.
 Then there exists a twisted form of $\mu_p$ over the function field $k(E)$ 
 that does not occur as subgroup scheme of an elliptic curve over $k(E)$.
 In particular, Question (A) has a negative answer for $k(E)$.
\end{Theorem}

\prf
We consider the morphism $E\to E$ which is given by multiplication 
with $p-1$.
This induces a field extension $k(E)\subset k(E)$, which is Galois with group
$(\ZZ/(p-1)\ZZ)^2$.
In particular, there exists a subfield $k(E)\subset L\subset k(E)$, such
that $L/k(E)$ is Galois with group $\ZZ/(p-1)\ZZ$.
Twisting $\mu_p$ with the $\mu_{p-1}$-torsor
$\Spec L\to\Spec k(E)$, we obtain a twisted form $\widetilde{\mu}_p$
of $\mu_p$.
We note that $L$ is the function field of an elliptic curve over $k$,
which is in particular a geometrically irreducible curve over $k$.

Assume there exists an elliptic curve $X$ over $k(E)$ containing
$\widetilde{\mu}_p$ as subgroup scheme.

As in the proof of Theorem \ref{main result} one first shows that
$j(X)\not\in k$.
Otherwise there exists a quadratic (resp. quartic, resp. sextic) 
twist $\widetilde{\mu}_p'$ of $\widetilde{\mu}_p$ which is defined over $k$,
which implies that the extension $L/k(E)$ is given by taking
the $(p-1)$.st root of an element of the form $a D^m$ with
$m=(p-1)/2$ (resp. $m=(p-1)/4$, resp. $m=(p-1)/6$) and $a\in k$.
However, this implies that $L$ is the function field of a curve
over $k$ which is not geometrically irreducible, a contradiction.

From $j(X)\not\in k$, we infer again that $\Spec L$ maps
dominantly to $\Ig(p)^{\rm ord}$, and using similar arguments as in the
proof of Theorem \ref{main result} we conclude 
that $\overline{\Ig(p)^{\rm ord}}$ has genus at most one,
which contradicts Proposition \ref{igusa proposition}.
\qed\medskip

If we choose $k$ to be a finite field in the previous 
two theorems, then we obtain examples of
global fields $K$ over which there exists
a twisted form $\widetilde{\mu}_p$
of $\mu_p$ that cannot be realised as subgroup scheme
of an elliptic curve over $K$.
On the other hand, Theorem \ref{local positive} tells us
that for every place $v$ of $K$ the group scheme
$\widetilde{\mu}_p\times_{\Spec K}\Spec K_v$ over the
completion $K_v$ can be realised as subgroup of an
elliptic curve over $K_v$.
Thus, although there exists a Hasse principle for
twisted forms over global fields by 
Proposition \ref{hasse principle}, we note:

\begin{Corollary}[No Hasse principle]
 \label{no hasse principle}
 Let $K$ be the function field of $\PP^1$ over a finite field of
 characteristic $p\geq13$ or the function field of an elliptic curve
 over a finite field of characteristic $p\geq17$.
 Then there is no Hasse principle for realising twisted forms
 of $\mu_p$ over $K$ as subgroup schemes of elliptic curves.
 \qed
\end{Corollary}

In view of the negative results Theorem \ref{main result} 
and Theorem \ref{main result elliptic} one can ask whether
one can realise more twisted forms of $\mu_p$ by also allowing
singular or nonproper curves carrying group structures rather
than only elliptic curves.
This turns out not to help:

If we do not insist on properness (but on geometric integrality), 
we are led to considering twisted forms of $\GG_m$.
If we do not insist on smoothness (but on geometric integrality)
then we are led to considering projective curves $C$ 
with $h^1(C,\OO_C)=1$ and having one singularity.
The smooth locus of $C$ is a twisted form of $\GG_m$ if
the singularity is a node and a twisted form of
$\GG_a$ if the singularity is a cusp.
Thus, we find twisted forms of $\mu_p$ also on nodal curves.

However, the group scheme of automorphisms of $\GG_m$
fixing the neutral element is $\ZZ/2\ZZ$ generated by
$t\mapsto t^{-1}$.
Thus, on twisted forms of nodal projective curves or on
twisted forms of $\GG_m$ we can only realise quadratic
twists of $\mu_p$.
Using quadratic twists of an elliptic curve $E$ 
containing $\mu_p$ as subgroup scheme, which always exists by
Remark \ref{mup ok}, we see that all quadratic twists of 
$\mu_p$ can be realised as subgroup schemes of elliptic
curves.
Thus we have shown

\begin{Proposition}
 \label{singular and nonproper}
 Let $k$ be field of positive characteristic $p$.
 Then every twisted form of $\mu_p$ that occurs on a twisted
 form of $\GG_m$ or on a twisted form of the singular nodal
 curve of arithmetic genus one can also be realised as
 subgroup scheme of an elliptic curve over $k$. \qed
\end{Proposition}

\section{The $p$-torsion subgroup scheme}
\label{sec:ptorsion}

We have seen in the previous sections that Question (A) has a positive
answer for local fields (Theorem \ref{local positive}), as well as
for all fields of characteristic $p\leq11$ (Corollary \ref{cor easy main result}).
We will see in this section that for those fields also Question (B) 
has a positive answer.

We need two technical lemmas to start with.

\begin{Lemma}
 \label{technical first}
 Let $G:=\mu_p\oplus(\ZZ/p\ZZ)$, $\Phi:G\times G\to\mu_p$ and 
 $\shA=\shAut(G,\Phi)$ as in Section \ref{sec:generalities}.
 Let $k$ be a field of positive characteristic $p$ and
 let $\Spec L\to\Spec k$ be an $\shA$-torsor.
 We denote by $(G,\Phi)\wedge^{\shA}\Spec L$ the twist of $(G,\Phi)$
 by this $\cal A$-torsor.
 Then
 \begin{enumerate}
  \item as a $k$-algebra $L$ is isomorphic to $k[x,y]/(x^p-a,y^{p-1}-b)$
    for some elements $a,b\in k$.
  \item as a scheme, the twist $(G,\Phi)\wedge^{\shA}\Spec L$
    is isomorphic to the spectrum of $k\oplus k[y]/(y^{p-1}-b)\oplus L$.
 \end{enumerate}
 Moreover, a $p^2$-dimensional $k$-algebra can carry at most one structure
 of a Hopf algebra making its spectrum into a twisted form of $(G,\Phi)$.
\end{Lemma}

\prf
We have seen in Section \ref{sec:generalities} that 
$\shA\iso\mu_{p-1}\rtimes\mu_p$.
Taking $\mu_{p-1}$-invariants we can factor $\Spec L\to\Spec k$
as $k\subset k':=k[x]/(x^p-a)$ for some $a\in k$ and
then $L=k'[z]/(z^{p-1}-c)$ for some $c\in k'$.
But then $b:=c^p$ lies in $k$ and thus
$k[x,y]/(x^p-a,y^{p-1}-b)$ is contained in $L$.
Comparing their dimensions as $k$-vector spaces they
must be equal and we get the first assertion.

The action of $\shA$ on $G$ is described in \cite[Section 2]{ls}.
This action has three orbits:
one orbit consists of $\Spec k$ (corresponding to the zero section),
one orbit has length $(p-1)$ (corresponding to $\ZZ/p\ZZ$ minus the
zero section) with infinitesimal isotropy groups, and finally one
orbit has length $p(p-1)$ upon which $\shA$ acts without fixed scheme.
From this we get the second assertion.

Finally, we note that every twisted form of $(G,\Phi)$ is autodual, i.e.,
isomorphic to its own Cartier dual.
Now, the coalgebra structure of a commutative and co-commutative Hopf algebra
is determined by its algebra structure, confer \cite[Section 2.4]{water}.
Hence a $k$-algebra can carry at most one structure of a Hopf algebra
making its spectrum into a twisted form of $(G,\Phi)$.
\qed

\begin{Lemma}
 \label{technical second}
 Let $E$ be an elliptic curve over a field of positive characteristic $p$
 with Hasse invariant $[h]\in k^\times/k^{\times(p-1)}$ and $j$-invariant
 $j(E)$.
 Then, as a scheme the $p$-torsion subgroup scheme $E[p]$ is isomorphic
 to the spectrum of 
 $k\oplus M\oplus L$, where $M=k[y]/(y^{p-1}-h)$ and $L=M[x]/(x^p-j(E))$.
\end{Lemma}

\prf
As a scheme, the twisted form $\ker(V)$ of $\ZZ/p\ZZ$, which is nothing
but the
\'etale quotient of $E[p]$ (compare with (\ref{connectedetale})),
is isomorphic to the spectrum of $k\oplus k[y]/(y^{p-1}-h)$.
These two summands corespond to the two $\shA$-orbits of
length $1$ and $p-1$ that we have seen in the proof
of Lemma \ref{technical first}.

To determine the third summand, we may assume $j(E)\not\in k^p$
for otherwise $E[p]$ is the direct sum of $\ker(V)$ and its
Cartier dual and the result is true in this case.
We note that $E[p]$
becomes isomorphic to $\mu_p\oplus(\ZZ/p\ZZ)$ over
$k[x,y]/(x^p-j(E),y^{p-1}-h)$.
From this it follows easily that $k[y]/(y^p-j(E))$
is contained in this third summand, which implies that
it contains $k[x,y]/(x^p-j(E),y^{p-1}-h)$.
Since this latter $k$-algebra is $p(p-1)$-dimensional 
like the summand we are looking for, they have to
coincide.
\qed\medskip

To answer Question (B) and to avoid trivialities, 
we may assume by Remark \ref{perfect remark}
that the field we are dealing with is not perfect.

\begin{Proposition}
 \label{hasse and j}
 Let $k$ be a nonperfect field of characteristic $p$.
 Then Question (B) has a positive answer for $k$ if for every 
 $[h]\in k^\times/k^{\times(p-1)}$ and every purely inseparable field extension
 $k\subset k'$ of degree $p$, there exists an elliptic curve $E$ over
 $k$ with Hasse invariant $[h]$ and $j$-invariant $j(E)$ such that
 $k'=k[\sqrt[p]{j(E)}]$.
\end{Proposition}

\prf
Let $\widetilde{G}$ be a twisted form of $(G,\Phi)$ over $k$.
By Lemma \ref{technical first} there exist 
$a,b\in k$ such that the scheme underlying $\widetilde{G}$
is isomorphic to the spectrum of $k\oplus M\oplus L$, where
$M=k[y]/(y^{p-1}-b)$ and $L=M[x]/(x^p-a)$.
By Lemma \ref{technical second} and our assumptions there exists
an elliptic curve $E$ over $k$
such that the scheme underlying $E[p]$ also 
is isomorphic to the spectrum of $k\oplus M\oplus L$.
But then, again by Lemma \ref{technical first}, 
the two schemes $\widetilde{G}$ and $E[p]$ are also
isomorphic as group schemes.
\qed

\begin{Theorem}
 Let $k$ be a field of positive characteristic $p\leq11$.
 Then every twisted form of $\mu_p\oplus(\ZZ/p\ZZ)$ endowed with a bilinear,
 alternating and nondegenerate pairing can be realised as subgroup scheme of
 an elliptic curve over $k$.
 In particular, Question (B) has a positive answer for $k$.
\end{Theorem}

\prf
By Corollary \ref{cor easy main result} and Remark \ref{perfect remark}
we may assume that $k$ is not perfect.
Given a purely inseparable field extension $k\subset k'$ of degree $p$,
there exists an element $j\in k$ such that $k'=k[\sqrt[p]{j}]$.
Then we choose an elliptic curve $E$ over $k$ with $j(E)=j$.

By Proposition \ref{hasse and j} it remains to realise all possible
Hasse invariants on such curves.
In characteristic $p=2$ there is only the Hasse invariant $h=[1]$ and
we get our result in this characteristic.
In characteristic $p=3$ we can use quadratic twists of this curve
$E$ to realise every possible Hasse invariant 
(see Remark \ref{refined easy main result})
and since twisting
does not change the $j$-invariant, we also get our result in
this characteristic.

In characteristic $p=5$ we may assume that elliptic curves are given
by a Weierstra\ss\ equation of the form $y^2=x^3+ax+b$.
The Hasse invariant of such a curve is $[2a]\in k^\times/k^{\times(p-1)}$ 
(see also the proof of Theorem \ref{easy main result}) and
we can realise every class of $k^\times/k^{\times(p-1)}$
by the Hasse invariant of an elliptic curve.
Replacing $a$ by $a^p$ does not change the Hasse invariant.
But then a straight forward computation shows that $k[\sqrt[p]{j(E)}]$ 
coincides with $k[\sqrt[p]{b}]$.
Choosing thus $a$ and $b$ appropriately, we see that the assumptions
of Proposition \ref{hasse and j} are fulfilled and our statement
follows for characteristic $5$.

The remaining cases $p=7$ and $p=11$ are similar to $p=5$ and therefore
left to the reader.
\qed

\begin{Theorem}
  Let $k$ be a perfect field of positive characteristic $p$.
 Then every twisted form of $\mu_p\oplus(\ZZ/p\ZZ)$ endowed with a bilinear,
 alternating and nondegenerate pairing can be realised as subgroup scheme of
 an elliptic curve over $k((t))$.
 In particular, Question (B) has a positive answer for $k((t))$.
\end{Theorem}

\prf
By Theorem \ref{local positive}, given a class in 
$k((t))^\times/k((t))^{\times(p-1)}$ there exists an elliptic curve $E$
over $k((t))$ having this class as Hasse invariant.
Moreover, looking at the proof of this result we see that we may assume
that  this elliptic curve has a $j$-invariant that does not lie in $k$.
If $j(E)$ lies in $k((t))^p$ then there exists an elliptic curve
$E'$ over $k((t))$ such that $E=E'^{(p)}$.
It is easy to see that the Hasse invariants of $E$ and $E'$
coincide.
We may thus assume that we have realised the given Hasse invariant
on an elliptic curve $E$ with $j(E)\not\in k((t))^p$.
Since there is only one inseparable extension $K'$ of $k((t))$ of
degree $p$ (here we use that $k$ is perfect), 
we must have $K'=k((t))[\sqrt[p]{j(E)}]$.
Applying Proposition \ref{hasse and j} our result follows.
\qed

\section{Supersingular elliptic curves}
\label{sec:supersingular}

Finally, we describe $p$-torsion subgroup schemes of supersingular elliptic curves, which turn
out to be much simpler than those of ordinary elliptic curves.

The kernel of Frobenius of a supersingular elliptic curve is a twisted form of $\alpha_p$.
Since $\shAut(\alpha_p)=\GG_a$ is a smooth group scheme, $\alpha_p$ does not possess
twisted forms over fields by Hilbert 90.
More generally, a twisted form of $\alpha_p$ over a scheme $S$ of positive
characteristic $p$ corresponds to a class in $H^1(S,\GG_a)$.
Since classes in this cohomology group can be trivialised Zariski-locally,
twisted forms of $\alpha_p$ over $S$ can be trivialised on a Zariski open
cover of $S$.

Over perfect fields, there exists only one non-split extension of $\alpha_p$
by itself that is autodual, namely $M_2$ in the notation of \cite[Section II.15.5]{oort}
\begin{equation}
 \label{m2sequence}
 0\,\to\,\alpha_p\,\to\,M_2\,\to\,\alpha_p\,\to\,0\,.
\end{equation}
Hence the $p$-torsion subgroup scheme of a supersingular elliptic curve is a twisted form
of $M_2$.
In particular, $M_2$ plays the role that $\mu_p\oplus(\ZZ/p\ZZ)$ plays for ordinary
elliptic curves.
Since $\alpha_p$ does not possess twisted forms over fields, twisted forms of $M_2$
correspond to twisted splittings of (\ref{m2sequence}), which are parametrised by 
$\Hfl{1}(\Spec k, \shHom(\alpha_p,\alpha_p))$, confer \cite[Chapter III.\S6.3.5]{dg}.
On the other hand, $\shHom(\alpha_p,\alpha_p)\iso\GG_a$ is a smooth group scheme, and 
hence this cohomology group is trivial for fields by Hilbert 90.
We deduce that $M_2$ does not possess twisted forms over fields.
Arguing as above, we conclude that a twisted form of $M_2$ over an arbitrary scheme $S$ of
positive characteristic $p$ can be trivialised on a Zariski-open cover of
$S$.

We have thus shown that the questions analogous to Questions (A) and (B)
trivially hold true for supersingular elliptic curves:

\begin{Theorem}
 Let $k$ be a field of positive characteristic $p$.
 \begin{enumerate}
  \item The kernel of Frobenius of a supersingular elliptic curve over $k$ is isomorphic
     to $\alpha_p$.
     The group scheme $\alpha_p$ does not possess twisted forms over $k$.
  \item The $p$-torsion subgroup scheme of a supersingular elliptic curve over $k$ is
     isomorphic to $M_2$.
     The group scheme $M_2$ does not possess twisted forms over $k$.
 \end{enumerate}
 More generally, every twisted form of $\alpha_p$, resp. of $M_2$, over an
 arbitrary scheme of positive characteristic $p$ can be embedded Zariski-locally
 into an elliptic curve.
\end{Theorem}

\end{document}